# Algebra of Fractions of Algebra with Conjugation

Aleks Kleyn

ABSTRACT. In the paper, I considered construction of algebra of fractions of algebra with conjugation. I also considered algebra of polynomials and algebra of rational mappings over algebra with conjugation.

## Contents

## 1. Auxiliary Theorems

**Theorem 1.1.** *Let* $a, b \in \operatorname{Im} A$. *Then*

$$\operatorname{Re}(ab) = \operatorname{Re}(ba) \tag{1.1}$$

$$\operatorname{Im}(ab) = -\operatorname{Im}(ba) \tag{1.2}$$

$$ab = (ba)^* \tag{1.3}$$

*Proof.* Since condition of lemma is true, then

$$a^* = -a \quad b^* = -b \tag{1.4}$$

The equation (1.3) follows from the equation (1.4) and the equation

$$ab = a^* b^* = (ba)^*$$

Therefore, equations (1.1), (1.2) follow from equations

$$ab = \operatorname{Re}(ab) + \operatorname{Im}(ab)$$
$$ab = (ba)^* = \operatorname{Re}(ba) - \operatorname{Im}(ba)$$

□

Aleks_Kleyn@MailAPS.org.
http://sites.google.com/site/AleksKleyn/.
http://arxiv.org/a/kleyn_a_1.
http://AleksKleyn.blogspot.com/.

**Theorem 1.2.** *Let $a \in A$. Then*

(1.5) $$aa^* = a^* a$$

*Proof.* The equation (1.5) follows from equations

$$aa^* = (\operatorname{Re} a)^2 - (\operatorname{Re} a)(\operatorname{Im} a) + (\operatorname{Im} a)(\operatorname{Re} a) - (\operatorname{Im} a)^2 = (\operatorname{Re} a)^2 - (\operatorname{Im} a)^2$$

$$a^* a = (\operatorname{Re} a)^2 + (\operatorname{Re} a)(\operatorname{Im} a) - (\operatorname{Im} a)(\operatorname{Re} a) - (\operatorname{Im} a)^2 = (\operatorname{Re} a)^2 - (\operatorname{Im} a)^2$$

□

**Theorem 1.3.** *Let $A$ be associative algebra with conjugation. Then*[1]

(1.6) $$(ab)(ab)^* = (aa^*)(bb^*)$$

*Proof.* Since

$$a = \operatorname{Re} a + \operatorname{Im} a$$
$$b = \operatorname{Re} b + \operatorname{Im} b$$

then

(1.7) $\ aa^* = (\operatorname{Re} a)^2 - (\operatorname{Re} a)(\operatorname{Im} a) + (\operatorname{Im} a)(\operatorname{Re} a) - (\operatorname{Im} a)^2 = (\operatorname{Re} a)^2 - (\operatorname{Im} a)^2$

(1.8) $\ bb^* = (\operatorname{Re} b)^2 - (\operatorname{Im} b)^2$

(1.9) $\quad ab = (\operatorname{Re} a)(\operatorname{Re} b) + (\operatorname{Re} a)(\operatorname{Im} b) + (\operatorname{Im} a)(\operatorname{Re} b) + (\operatorname{Im} a)(\operatorname{Im} b)$

$\qquad (ab)^* = (\operatorname{Re} a)(\operatorname{Re} b) - (\operatorname{Re} a)(\operatorname{Im} b) - (\operatorname{Im} a)(\operatorname{Re} b) + ((\operatorname{Im} a)(\operatorname{Im} b))^*$

(1.10) $\qquad\quad = (\operatorname{Re} a)(\operatorname{Re} b) - (\operatorname{Re} a)(\operatorname{Im} b) - (\operatorname{Im} a)(\operatorname{Re} b) + (\operatorname{Im} b)^*(\operatorname{Im} a)^*$

$\qquad\quad = (\operatorname{Re} a)(\operatorname{Re} b) - (\operatorname{Re} a)(\operatorname{Im} b) - (\operatorname{Im} a)(\operatorname{Re} b) + (\operatorname{Im} b)(\operatorname{Im} a)$

$\qquad\quad (ab)(ab)^*$

$\qquad\quad = ((\operatorname{Re} a)(\operatorname{Re} b) + (\operatorname{Re} a)(\operatorname{Im} b) + (\operatorname{Im} a)(\operatorname{Re} b) + (\operatorname{Im} a)(\operatorname{Im} b))$

$\qquad\qquad * ((\operatorname{Re} a)(\operatorname{Re} b) - (\operatorname{Re} a)(\operatorname{Im} b) - (\operatorname{Im} a)(\operatorname{Re} b) + (\operatorname{Im} b)(\operatorname{Im} a))$

$\qquad\quad = (\operatorname{Re} a)(\operatorname{Re} b)(\operatorname{Re} a)(\operatorname{Re} b) + (\operatorname{Re} a)(\operatorname{Im} b)(\operatorname{Re} a)(\operatorname{Re} b)_{-1-}$

$\qquad\qquad + (\operatorname{Im} a)(\operatorname{Re} b)(\operatorname{Re} a)(\operatorname{Re} b)_{-2-} + (\operatorname{Im} a)(\operatorname{Im} b)(\operatorname{Re} a)(\operatorname{Re} b)_{-3-}$

(1.11) $\qquad\quad - (\operatorname{Re} a)(\operatorname{Re} b)(\operatorname{Re} a)(\operatorname{Im} b)_{-1-} - (\operatorname{Re} a)(\operatorname{Im} b)(\operatorname{Re} a)(\operatorname{Im} b)$

$\qquad\qquad - (\operatorname{Im} a)(\operatorname{Re} b)(\operatorname{Re} a)(\operatorname{Im} b)_{-3-} - ((\operatorname{Im} a)(\operatorname{Im} b))(\operatorname{Re} a)(\operatorname{Im} b)$

$\qquad\qquad - (\operatorname{Re} a)(\operatorname{Re} b)(\operatorname{Im} a)(\operatorname{Re} b)_{-2-} - (\operatorname{Re} a)(\operatorname{Im} b)(\operatorname{Im} a)(\operatorname{Re} b)_{-4-}$

$\qquad\qquad - (\operatorname{Im} a)(\operatorname{Re} b)(\operatorname{Im} a)(\operatorname{Re} b) - ((\operatorname{Im} a)(\operatorname{Im} b))(\operatorname{Im} a)(\operatorname{Re} b)$

$\qquad\qquad + (\operatorname{Re} a)(\operatorname{Re} b)(\operatorname{Im} b)(\operatorname{Im} a)_{-4-} + (\operatorname{Re} a)(\operatorname{Im} b)((\operatorname{Im} b)(\operatorname{Im} a))$

$\qquad\qquad + (\operatorname{Im} a)(\operatorname{Re} b)((\operatorname{Im} b)(\operatorname{Im} a)) + ((\operatorname{Im} a)(\operatorname{Im} b))((\operatorname{Im} b)(\operatorname{Im} a))$

---

[1]The theorem 1.3 has more simple prove. Namely, the theorem follows from the equation
$$(ab)(ab)^* = (ab)(b^* a^*) = a((bb^*)a^*) = (aa^*)(bb^*)$$
However, I hope to find conditions when the theorem 1.3 is true for non associative algebra. The proof in the text is the basis for future research.

From the equation (1.11), it follows that

$$
\begin{aligned}
(1.12)\quad (ab)(ab)^* &= (\operatorname{Re} a)^2 (\operatorname{Re} b)^2 \\
&\quad - (\operatorname{Re} a)^2 (\operatorname{Im} b)^2 - (\operatorname{Im} a)^2 (\operatorname{Re} b)^2 \\
&\quad + ((\operatorname{Im} a)(\operatorname{Im} b))((\operatorname{Im} b)(\operatorname{Im} a)) \\
&\quad - ((\operatorname{Im} a)(\operatorname{Im} b))(\operatorname{Re} a)(\operatorname{Im} b) + (\operatorname{Re} a)(\operatorname{Im} b)((\operatorname{Im} b)(\operatorname{Im} a)) \\
&\quad - ((\operatorname{Im} a)(\operatorname{Im} b))(\operatorname{Im} a)(\operatorname{Re} b) + (\operatorname{Im} a)(\operatorname{Re} b)((\operatorname{Im} b)(\operatorname{Im} a))
\end{aligned}
$$

Since the algebra $A$ is associative, then

$$
\begin{aligned}
(1.13)\quad ((\operatorname{Im} a)(\operatorname{Im} b))(\operatorname{Im} b) &= (\operatorname{Im} a)((\operatorname{Im} b)(\operatorname{Im} b)) \\
&= ((\operatorname{Im} b)(\operatorname{Im} b))(\operatorname{Im} a) \\
&= (\operatorname{Im} b)((\operatorname{Im} b)(\operatorname{Im} a))
\end{aligned}
$$

$$
(1.14)\quad ((\operatorname{Im} a)(\operatorname{Im} b))(\operatorname{Im} a) = (\operatorname{Im} a)((\operatorname{Im} b)(\operatorname{Im} a))
$$

$$
\begin{aligned}
(1.15)\quad ((\operatorname{Im} a)(\operatorname{Im} b))((\operatorname{Im} b)(\operatorname{Im} a)) &= (\operatorname{Im} a)((\operatorname{Im} b)((\operatorname{Im} b)(\operatorname{Im} a))) \\
&= (\operatorname{Im} a)(((\operatorname{Im} b)(\operatorname{Im} b))(\operatorname{Im} a)) \\
&= ((\operatorname{Im} a)(\operatorname{Im} a))((\operatorname{Im} b)(\operatorname{Im} b)) \\
&= (\operatorname{Im} a)^2 (\operatorname{Im} b)^2
\end{aligned}
$$

It follows from equations (1.12), (1.13), (1.14), (1.15) that

$$
\begin{aligned}
(1.16)\quad (ab)(ab)^* &= (\operatorname{Re} a)^2 (\operatorname{Re} b)^2 - (\operatorname{Re} a)^2 (\operatorname{Im} b)^2 \\
&\quad - (\operatorname{Im} a)^2 (\operatorname{Re} b)^2 + (\operatorname{Im} a)^2 (\operatorname{Im} b)^2 \\
&= (\operatorname{Re} a)^2 ((\operatorname{Re} b)^2 - (\operatorname{Im} b)^2) \\
&\quad - (\operatorname{Im} a)^2 ((\operatorname{Re} b)^2 - (\operatorname{Im} b)^2) \\
&= ((\operatorname{Re} a)^2 - (\operatorname{Im} a)^2)((\operatorname{Re} b)^2 - (\operatorname{Im} b)^2)
\end{aligned}
$$

The equation (1.6) follows from equations (1.7), (1.8), (1.16). □

**Theorem 1.4.** *Let $A$ be associative algebra with conjugation. Then*

$$
(1.17)\quad \left(\prod_{i=1}^m a_i\right)\left(\prod_{i=1}^m a_i\right)^* = \prod_{i=1}^m (a_i a_i^*)
$$

*Proof.* For $m=1$, the theorem is obvious. For $m=2$, the theorem follows from the theorem 1.3. Let the theorem is true for $m = p-1$. Let

$$
b = \prod_{i=1}^{p-1} a_i
$$

Then

$$\left(\prod_{i=1}^{p} a_i\right)\left(\prod_{i=1}^{p} a_i\right)^* = \left(\left(\prod_{i=1}^{p-1} a_i\right) a_p\right)\left(\left(\prod_{i=1}^{p-1} a_i\right) a_p\right)^*$$
$$= (ba_p)(ba_p)^* = (ba_p)(a_p^* b^*)$$
$$= b(a_p a_p^*) b^* = bb^*(a_p a_p^*)$$
$$= \left(\prod_{i=1}^{p-1} a_i\right)\left(\prod_{i=1}^{p-1} a_i\right)^* (a_p a_p^*)$$
$$= \left(\prod_{i=1}^{p-1} (a_i a_i^*)\right)(a_p a_p^*) = \prod_{i=1}^{p} (a_i a_i^*)$$

Therefore, the theorem is true for $m = p$. □

## 2. Field of Fractions of Scalar Algebra

Let $D$ be commutative ring. Let $A$ be $D$-algebra with conjugation. According to the definition [7]-4.2, scalar algebra $\operatorname{Re} A$ is commutative associative ring. Let the ring $\operatorname{Re} A$ be entire.[2] Then there exists field $F$ of fractions of ring $\operatorname{Re} A$.[3]

According to construction that was done in subsections [8]-4.4.2, [8]-4.4.3, a diagram of representations of $\operatorname{Re} A$-algebra $A$ has form

$$\operatorname{Re} A \xrightarrow{f_{1,2}}_{*} A \xrightarrow{f_{2,3}}_{*} A \qquad f_{1,2}(d) : a \to d\,a$$
$$\uparrow_{*f_{1,2}} \qquad f_{2,3}(a) : b \to C_A(a,b)$$
$$\operatorname{Re} A \qquad C_A \in \mathcal{L}(A^2; A)$$

A diagram of representations of $F$-algebra $G$ has form

$$F \xrightarrow{g_{1,2}}_{*} B \xrightarrow{g_{2,3}}_{*} B \qquad g_{1,2}(d) : a \to d\,a$$
$$\uparrow_{*g_{1,2}} \qquad g_{2,3}(a) : b \to C_B(a,b)$$
$$F \qquad C_B \in \mathcal{L}(B^2; B)$$

We define $F$-algebra $B$ such that there exists linear homomorphism[4] of $\operatorname{Re} A$-algebra $A$ into $F$-algebra $B$

$$r_1 : \operatorname{Re} A \to F \quad r_2 : A \to B$$

such that ring homomorphism $r_1$ is embedding of the ring $\operatorname{Re} A$ into the ring $F$ (page [1]-108)

$$(2.1) \qquad r_1(d) = d/1$$

and image of a basis $\overline{\overline{e}}_A$ of $\operatorname{Re} A$-module $A$ under mapping $r_2$ is a basis of $F$-vector space $B$

$$(2.2) \qquad r_2 \circ \overline{e}_{A \cdot i} = \overline{e}_{B \cdot i}$$

---

[2] See the definition of entire ring on the page [1]-91.
[3] Construction of field of fractions is considered in [1], pages 107 - 110.
[4] See the definition [6]-6.2.

Based on the equation (2.1) we can identify $d \in \operatorname{Re} A$ and its image

$$r_1(d) = d$$

**Theorem 2.1.** *$F$-algebra $B$ is algebra with conjugation. The field $F$ is scalar algebra of $F$-algebra $B$. Structural constants of $F$-algebra $B$ coincide with structural constants of $\operatorname{Re} A$-algebra $A$*

$$(2.3) \qquad C_{B \cdot ij}^{\phantom{B\cdot}k} = C_{A \cdot ij}^{\phantom{A\cdot}k}$$

*Proof.* From the equation (2.2) it follows that

$$(2.4) \qquad r_{2 \cdot i}^{\phantom{2\cdot}k} = \delta_i^k$$

From the equation (2.4) and the theorem [6]-6.4 it follows that

$$(2.5) \qquad C_{A \cdot ij}^{\phantom{A\cdot}l} = r_1(C_{A \cdot ij}^{\phantom{A\cdot}k})\delta_k^l = \delta_i^p \delta_j^q C_{B \cdot pq}^{\phantom{B\cdot}l}$$

The equation (2.2) follows from the equation (2.5). From the equation (2.5) and the theorem [7]-4.5 it follows that $F$-algebra $B$ is algebra with conjugation and the field $F$ is scalar algebra of algebra $B$. $\square$

Below we will assume that $\operatorname{Re} A$ is a field.

According to the theorem [7]-4.9

$$(2.6) \qquad aa^* \in \operatorname{Re} A$$

In contrast to complex field and quaternion algebra, the field $\operatorname{Re} A$ can be different from the real field. The concept of order may be missing in the field $\operatorname{Re} A$. So we cannot accept the expression (2.6) as norm in the algebra $A$. Even more, this expression can be equal 0.

**Theorem 2.2.** *$a \in A$ is invertible in the $\operatorname{Re} A$-algebra $A$ iff*

$$(2.7) \qquad aa^* \neq 0$$

*Since the condition (2.6) is true, then*

$$(2.8) \qquad a^{-1} = \frac{1}{aa^*} a^*$$

*Proof.* By definition, $a \in A$ is invertible, if there exists $a^{-1} \in A$ such that

$$(2.9) \qquad aa^{-1} = 1$$

From the equation

$$aa^* = aa^*$$

and the equation (2.6), it follows that

$$(2.10) \qquad \frac{1}{aa^*} aa^* = \frac{1}{aa^*} aa^*$$

We can represent right part of the equation (2.10) as

$$(2.11) \qquad \frac{1}{aa^*} aa^* = \frac{1}{aa^*} \frac{aa^*}{1} = 1$$

Since the product in $\operatorname{Re} A$-algebra $A$ is bilinear mapping, then we can represent left part of the equation (2.10) as

$$(2.12) \qquad \frac{1}{aa^*} aa^* = a \left( \frac{1}{aa^*} a^* \right)$$

From equations (2.10), (2.11), (2.12), it follows that

$$a\left(\frac{1}{aa^*}a^*\right) = 1 \tag{2.13}$$

(2.8) follows from equations (2.9), (2.13).

Since

$$aa^* = 0 \tag{2.14}$$

then $a \in A$ is left and right zero divisor.[5] According to the theorem [3]-6.3, $a \in A$ does not have inverse. □

*Remark* 2.3. Using the notation considered in the beginning of this section we can say that it follows from the theorem 2.2 that $F$-algebra $B$ is algebra which has the greatest possible set of invertible elements of Re $A$-algebra $A$. $F$-algebra $B$ is called **algebra of fractions of algebra with conjugation** $A$. □

**Definition 2.4.** Let us denote

$$A_0 = \{a \in A : a\,a^* = 0\}$$

**the set of zeros of algebra** $A$. According to the theorem 2.2, $a \in A_0$ iff either $a = 0$, or $a$ is zero divisor. Let us denote

$$A_1 = \{a \in A : a\,a^* \neq 0\}$$

**set of invertible elements of algebra** $A$. □

Let $a \in A_1$, $b \in A$. **Left fraction** is represented by expression

$$a^{-1}b = \frac{1}{a\,a^*}a^*b$$

**Right fraction** is represented by expression

$$ba^{-1} = \frac{1}{a\,a^*}b\,a^*$$

The set of fractions in algebra is not limited by left or right fractions. For instance, expressions

$$(a^{-1}b)(c^{-1}d) \quad a^{-1}b + c^{-1}d$$

are also fractions.

We can define few equivalence relations on the set of fractions. For instance, since $d \in \text{Re}\,A$, then

$$a^{-1}b = (da)^{-1}(db) \quad (a^{-1}(db))(c^{-1}f) = (a^{-1}b)(c^{-1}(df))$$

However, the question about canonical form of fraction is not trivial task, at least, at current time.

---

[5]From equations

$$aa^* = C^0_{00}a^0a^0 - C^0_{kl}a^k a^l$$
$$a^*a = C^0_{00}a^0a^0 - C^0_{kl}a^k a^l$$

and the equation (2.14), it follows that $aa^* = 0$.

## 3. Algebra of Polynomials

Let $D$ be the commutative ring of characteristic $0$. Let $A$ be $D$-algebra. **Algebra of polynomials** $A[x]$ **over** $D$-**algebra** $A$ is generated by the set of monomials. The following theorem (the section [5]-5.2) describes the structure of the monomial $p_k$ of power $k$, $k > 0$, in one variable over associative $D$-algebra $A$.

**Theorem 3.1.** *Monomial of power $0$ has form $a_0$, $a_0 \in A$. For $k > 0$,*

$$p_k(x) = p_{k-1}(x) x a_k$$

*where $a_k \in A$.*

*Proof.* Actually, last factor of monomial $p_k(x)$ is either $a_k \in A$, or has form $x^l$, $l \geq 1$. In the later case we assume $a_k = 1$. Factor preceding $a_k$ has form $x^l$, $l \geq 1$. We can represent this factor as $x^{l-1}x$. Therefore, we proved the statement. $\square$

In particular, monomial of power $1$ has form $p_1(x) = a_0 x a_1$. From theorems 3.1, [5]-3.41, it follows that we can associate the tensor $a_0 \otimes a_1 \otimes ... \otimes a_k$ to each monomial $p_k$.

Order of the factors is essential in the nonassociative algebra. So the theorem 3.1 gets following form.

**Theorem 3.2.** *Monomial of power $0$ has form $a_0$, $a_0 \in A$. For $k > 0$, there exist monomials $p_l$, $p_m$, $m + l = k$, such that*

$$p_k(x) = a_{k \cdot 1} p_l(x) a_{k \cdot 2} p_m(x) a_{k \cdot 3}$$

*where $a_{k \cdot 1}$, $a_{k \cdot 2}$, $a_{k \cdot 3} \in A$.* $\square$

Since $D$-algebra $A$ is algebra with conjugation, then we can extend the mapping of conjugation onto algebra of polynomials as well we can consider polynomials over ring $\operatorname{Re} A$. Since the structure of polynomial over ring $\operatorname{Re} A$ is different from the structure of polynomial over algebra $A$, then determination of relationship between algebras $\operatorname{Re} A[x]$ and $A[x]$ is nontrivial problem.

Since $p(x)$ is monomial over algebra, then, according to the theorem 1.4, we can consider an expression $p(x)(p(x))^*$ as polynomial $r(y)$ with variable $y = xx^*$ over algebra $\operatorname{Re} A$. Although for an arbitrary polynomial $p(x) \in A[x]$, expressions $p(x) + (p(x))^*$, $p(x)(p(x))^*$ take values in ring $\operatorname{Re} A$, it is not clear whether we can consider these expressions as polynomials over ring $\operatorname{Re} A$.

There exist algebras with conjugation where conjugation does not depend linearly on identity mapping (see, for instance, section [2]-6). In such case for any polynomial, expression

$$p(x)(p(x))^*$$

depends from two variables: $x$ and $x^*$.

Consider algebras with conjugation where conjugation linearly depends on identity mapping

$$x^* = s(x) = s_{i \cdot 0} \, x \, s_{i \cdot 1}$$

(for instance, the mapping [4]-(4.3.35) in quaternion algebra, the mapping [4]-(4.5.99) in octonion algebra). In such case for any polynomial $p(x)$, expression

$$p(x)(p(x))^* = r(x)$$

is polynomial.

Let the ring $\operatorname{Re} A$ be a field.

**Definition 3.3.** Let $p(x)$ be polynomial over algebra $A$. $a \in A$ is called **root of the polynomial** $p$, if $p(a) \in A_0$. □

According to the theorem 2.2, for polynomial $p(x)$, following expression is defined

$$\frac{1}{p(x)(p(x))^*} \in \operatorname{Re} A \tag{3.1}$$

for any $x \in A$ which is different from root of polynomial $p(x)$. Therefore, the mapping

$$(p(x))^{-1} = \frac{1}{p(x)(p(x))^*}(p(x))^* \tag{3.2}$$

is defined properly for $x$ which is not root of polynomial $p$.

Algebra $A(x)$ generated by expressions like (3.2) is called **algebra of rational mappings of algebra** $A$.

## 4. Ideal

**Definition 4.1.** Subgroup $B$ of additive group of algebra with conjugation $A$ is called **left ideal of algebra**,[6] if

$$aB \subset B \quad a \in A$$

Subgroup $B$ of additive group of algebra with conjugation $A$ is called **right ideal of algebra**, if

$$Ba \subset B \quad a \in A$$

Subgroup $B$ of additive group of algebra with conjugation $A$ is called **ideal of algebra**, if $B$ is both a left and a right ideal. □

**Example 4.2.** Let $A$ be algebra with conjugation, $a \in A$. The set $Aa$ is left ideal called **left principal ideal** of algebra $A$.[7] The set $aA$ is right ideal called **right principal ideal** of algebra $A$. The set $AaA$ is ideal called **principal ideal** of algebra $A$. □

**Theorem 4.3.** *Let $A$ be associative algebra with conjugation, $a \in A_0$. Then*

$$Aa \in A_0 \tag{4.1}$$

$$aA \in A_0 \tag{4.2}$$

*Proof.* Let $b \in A$. From the definition 2.4 it follows that

$$(ba)(ba)^* = (ba)(a^*b^*) = b(aa^*)b^* = 0 \tag{4.3}$$

$$(ab)(ab)^* = (ab)(b^*a^*) = a(bb^*)a^* = (bb^*)(aa^*) = 0 \tag{4.4}$$

The statement (4.1) follows from the equation (4.3). The statement (4.2) follows from the equation (4.4). □

---

[6]This definition is based on the definition [1], page 86.
[7][1], page 86.

**Theorem 4.4.** *Let $A$ be associative algebra with conjugation, $a \in A$, $b \in A_0$. Algebra of polynomials $A[x]$ has left ideal*

$$Z_l^1(a,b)A[x] = \{p \in A[x] : p(a) \in Ab\}$$

*Algebra of polynomials $A[x]$ has right ideal*

$$Z_r^1(a,b)A[x] = \{p \in A[x] : p(a) \in bA\}$$

*Proof.* The theorem follows from definitions considered in the example 4.2 and the theorem 4.3. □

We can prove similar theorems.

**Theorem 4.5.** *Let $A$ be associative algebra with conjugation, $a \in A$. Algebra of polynomials $A[x]$ has ideal*

$$Z^1(a)A[x] = \{p \in A[x] : p(a) = 0\}$$

□

**Theorem 4.6.** *Let $A$ be associative algebra with conjugation, $a \in A$, $b \in A_0$. Algebra of rational mappings $A(x)$ has left ideal*

$$Z_l^1(a,b)A(x) = \{p \in A(x) : p(a) \in Ab\}$$

*Algebra of rational mappings $A(x)$ has right ideal*

$$Z_r^1(a,b)A(x) = \{p \in A(x) : p(a) \in bA\}$$

□

**Theorem 4.7.** *Let $A$ be associative algebra with conjugation, $a \in A$. Algebra of rational mappings $A(x)$ has ideal*

$$Z^1(a)A(x) = \{p \in A(x) : p(a) = 0\}$$

□

## 6. Index

# 7. Special Symbols and Notations

# Алгебра частных алгебры с сопряжением

Александр Клейн

Аннотация. В статье рассмотренно возможность построения алгебры частных алгебры с сопряжением. Я также рассмотрел алгебру многочленов и алгебру рациональных отображений над алгеброй с сопряжением.

Содержание

## 1. Вспомогательные теоремы

**Теорема 1.1.** *Пусть $a$, $b \in \operatorname{Im} A$. Тогда*

$$(1.1) \qquad \operatorname{Re}(ab) = \operatorname{Re}(ba)$$

$$(1.2) \qquad \operatorname{Im}(ab) = -\operatorname{Im}(ba)$$

$$(1.3) \qquad ab = (ba)^*$$

*Доказательство.* Если условие леммы выполнено, то

$$(1.4) \qquad a^* = -a \quad b^* = -b$$

Равенство (1.3) является следствием равенства (1.4) и равенства

$$ab = a^*b^* = (ba)^*$$

Следовательно, равенства (1.1), (1.2) являются следствием равенств

$$ab = \operatorname{Re}(ab) + \operatorname{Im}(ab)$$
$$ab = (ba)^* = \operatorname{Re}(ba) - \operatorname{Im}(ba)$$

□

Aleks_Kleyn@MailAPS.org.
http://sites.google.com/site/AleksKleyn/.
http://arxiv.org/a/kleyn_a_1.
http://AleksKleyn.blogspot.com/.

**Теорема 1.2.** *Пусть $a \in A$. Тогда*

(1.5) $$aa^* = a^*a$$

*Доказательство.* Равенство (1.5) следует из равенств

$$aa^* = (\operatorname{Re} a)^2 - (\operatorname{Re} a)(\operatorname{Im} a) + (\operatorname{Im} a)(\operatorname{Re} a) - (\operatorname{Im} a)^2 = (\operatorname{Re} a)^2 - (\operatorname{Im} a)^2$$

$$a^*a = (\operatorname{Re} a)^2 + (\operatorname{Re} a)(\operatorname{Im} a) - (\operatorname{Im} a)(\operatorname{Re} a) - (\operatorname{Im} a)^2 = (\operatorname{Re} a)^2 - (\operatorname{Im} a)^2$$

$\square$

**Теорема 1.3.** *Пусть $A$ - ассоциативная алгебра с сопряжением. Тогда*[1]

(1.6) $$(ab)(ab)^* = (aa^*)(bb^*)$$

*Доказательство.* Поскольку

$$a = \operatorname{Re} a + \operatorname{Im} a$$
$$b = \operatorname{Re} b + \operatorname{Im} b$$

то

(1.7) $\quad aa^* = (\operatorname{Re} a)^2 - (\operatorname{Re} a)(\operatorname{Im} a) + (\operatorname{Im} a)(\operatorname{Re} a) - (\operatorname{Im} a)^2 = (\operatorname{Re} a)^2 - (\operatorname{Im} a)^2$

(1.8) $\quad bb^* = (\operatorname{Re} b)^2 - (\operatorname{Im} b)^2$

(1.9) $\quad ab = (\operatorname{Re} a)(\operatorname{Re} b) + (\operatorname{Re} a)(\operatorname{Im} b) + (\operatorname{Im} a)(\operatorname{Re} b) + (\operatorname{Im} a)(\operatorname{Im} b)$

$\quad\quad (ab)^* = (\operatorname{Re} a)(\operatorname{Re} b) - (\operatorname{Re} a)(\operatorname{Im} b) - (\operatorname{Im} a)(\operatorname{Re} b) + ((\operatorname{Im} a)(\operatorname{Im} b))^*$

(1.10) $\quad\quad\quad = (\operatorname{Re} a)(\operatorname{Re} b) - (\operatorname{Re} a)(\operatorname{Im} b) - (\operatorname{Im} a)(\operatorname{Re} b) + (\operatorname{Im} b)^*(\operatorname{Im} a)^*$

$\quad\quad\quad = (\operatorname{Re} a)(\operatorname{Re} b) - (\operatorname{Re} a)(\operatorname{Im} b) - (\operatorname{Im} a)(\operatorname{Re} b) + (\operatorname{Im} b)(\operatorname{Im} a)$

$\quad\quad (ab)(ab)^*$

$\quad\quad = ((\operatorname{Re} a)(\operatorname{Re} b) + (\operatorname{Re} a)(\operatorname{Im} b) + (\operatorname{Im} a)(\operatorname{Re} b) + (\operatorname{Im} a)(\operatorname{Im} b))$

$\quad\quad * ((\operatorname{Re} a)(\operatorname{Re} b) - (\operatorname{Re} a)(\operatorname{Im} b) - (\operatorname{Im} a)(\operatorname{Re} b) + (\operatorname{Im} b)(\operatorname{Im} a))$

$\quad\quad = (\operatorname{Re} a)(\operatorname{Re} b)(\operatorname{Re} a)(\operatorname{Re} b) + (\operatorname{Re} a)(\operatorname{Im} b)(\operatorname{Re} a)(\operatorname{Re} b)_{-1-}$

$\quad\quad + (\operatorname{Im} a)(\operatorname{Re} b)(\operatorname{Re} a)(\operatorname{Re} b)_{-2-} + (\operatorname{Im} a)(\operatorname{Im} b)(\operatorname{Re} a)(\operatorname{Re} b)_{-3-}$

(1.11) $\quad\quad - (\operatorname{Re} a)(\operatorname{Re} b)(\operatorname{Re} a)(\operatorname{Im} b)_{-1-} - (\operatorname{Re} a)(\operatorname{Im} b)(\operatorname{Re} a)(\operatorname{Im} b)$

$\quad\quad - (\operatorname{Im} a)(\operatorname{Re} b)(\operatorname{Re} a)(\operatorname{Im} b)_{-3-} - ((\operatorname{Im} a)(\operatorname{Im} b))(\operatorname{Re} a)(\operatorname{Im} b)$

$\quad\quad - (\operatorname{Re} a)(\operatorname{Re} b)(\operatorname{Im} a)(\operatorname{Re} b)_{-2-} - (\operatorname{Re} a)(\operatorname{Im} b)(\operatorname{Im} a)(\operatorname{Re} b)_{-4-}$

$\quad\quad - (\operatorname{Im} a)(\operatorname{Re} b)(\operatorname{Im} a)(\operatorname{Re} b) - ((\operatorname{Im} a)(\operatorname{Im} b))(\operatorname{Im} a)(\operatorname{Re} b)$

$\quad\quad + (\operatorname{Re} a)(\operatorname{Re} b)(\operatorname{Im} b)(\operatorname{Im} a)_{-4-} + (\operatorname{Re} a)(\operatorname{Im} b)((\operatorname{Im} b)(\operatorname{Im} a))$

$\quad\quad + (\operatorname{Im} a)(\operatorname{Re} b)((\operatorname{Im} b)(\operatorname{Im} a)) + ((\operatorname{Im} a)(\operatorname{Im} b))((\operatorname{Im} b)(\operatorname{Im} a))$

---

[1]Теорема 1.3 имеет более простое доказательство. А именно, теорема следует из равенства

$$(ab)(ab)^* = (ab)(b^*a^*) = a((bb^*)a^*) = (aa^*)(bb^*)$$

Однако, я надеюсь найти условия, когда теорема 1.3 верна для неассоциативной алгебры. Приведенное доказательство является основой для будущего исследования.

Из равенства (1.11) следует

$$
\begin{aligned}
&(ab)(ab)^* \\
=&(\operatorname{Re} a)^2(\operatorname{Re} b)^2 \\
&-(\operatorname{Re} a)^2(\operatorname{Im} b)^2 - (\operatorname{Im} a)^2(\operatorname{Re} b)^2 \\
&+((\operatorname{Im} a)(\operatorname{Im} b))((\operatorname{Im} b)(\operatorname{Im} a)) \\
&-((\operatorname{Im} a)(\operatorname{Im} b))(\operatorname{Re} a)(\operatorname{Im} b) + (\operatorname{Re} a)(\operatorname{Im} b)((\operatorname{Im} b)(\operatorname{Im} a)) \\
&-((\operatorname{Im} a)(\operatorname{Im} b))(\operatorname{Im} a)(\operatorname{Re} b) + (\operatorname{Im} a)(\operatorname{Re} b)((\operatorname{Im} b)(\operatorname{Im} a))
\end{aligned}
\qquad (1.12)
$$

Если произведение в алгебре $A$ ассоциативно, то

$$
\begin{aligned}
((\operatorname{Im} a)(\operatorname{Im} b))(\operatorname{Im} b) &= (\operatorname{Im} a)((\operatorname{Im} b)(\operatorname{Im} b)) \\
&= ((\operatorname{Im} b)(\operatorname{Im} b))(\operatorname{Im} a) \\
&= (\operatorname{Im} b)((\operatorname{Im} b)(\operatorname{Im} a))
\end{aligned}
\qquad (1.13)
$$

$$
((\operatorname{Im} a)(\operatorname{Im} b))(\operatorname{Im} a) = (\operatorname{Im} a)((\operatorname{Im} b)(\operatorname{Im} a)) \qquad (1.14)
$$

$$
\begin{aligned}
((\operatorname{Im} a)(\operatorname{Im} b))((\operatorname{Im} b)(\operatorname{Im} a)) &= (\operatorname{Im} a)((\operatorname{Im} b)((\operatorname{Im} b)(\operatorname{Im} a))) \\
&= (\operatorname{Im} a)(((\operatorname{Im} b)(\operatorname{Im} b))(\operatorname{Im} a)) \\
&= ((\operatorname{Im} a)(\operatorname{Im} a))((\operatorname{Im} b)(\operatorname{Im} b)) \\
&= (\operatorname{Im} a)^2(\operatorname{Im} b)^2
\end{aligned}
\qquad (1.15)
$$

Из равенств (1.12), (1.13), (1.14), (1.15) следует

$$
\begin{aligned}
(ab)(ab)^* =&(\operatorname{Re} a)^2(\operatorname{Re} b)^2 - (\operatorname{Re} a)^2(\operatorname{Im} b)^2 \\
&-(\operatorname{Im} a)^2(\operatorname{Re} b)^2 + (\operatorname{Im} a)^2(\operatorname{Im} b)^2 \\
=&(\operatorname{Re} a)^2((\operatorname{Re} b)^2 - (\operatorname{Im} b)^2) \\
&-(\operatorname{Im} a)^2((\operatorname{Re} b)^2 - (\operatorname{Im} b)^2) \\
=&((\operatorname{Re} a)^2 - (\operatorname{Im} a)^2)((\operatorname{Re} b)^2 - (\operatorname{Im} b)^2)
\end{aligned}
\qquad (1.16)
$$

Равенство (1.6) следует из равенств (1.7), (1.8), (1.16). □

**Теорема 1.4.** *Пусть $A$ - ассоциативная алгебра с сопряжением. Тогда*

$$
\left(\prod_{i=1}^{m} a_i\right)\left(\prod_{i=1}^{m} a_i\right)^* = \prod_{i=1}^{m}(a_i a_i^*) \qquad (1.17)
$$

*Доказательство.* Для $m=1$ утверждение теоремы очевидно. Для $m=2$ теорема следует из теоремы 1.3. Пусть утверждение теоремы верно для $m=p-1$. Пусть

$$
b = \prod_{i=1}^{p-1} a_i
$$

Тогда

$$\left(\prod_{i=1}^{p} a_i\right)\left(\prod_{i=1}^{p} a_i\right)^* = \left(\left(\prod_{i=1}^{p-1} a_i\right) a_p\right)\left(\left(\prod_{i=1}^{p-1} a_i\right) a_p\right)^*$$
$$= (ba_p)(ba_p)^* = (ba_p)(a_p^* b^*)$$
$$= b(a_p a_p^*) b^* = bb^*(a_p a_p^*)$$
$$= \left(\prod_{i=1}^{p-1} a_i\right)\left(\prod_{i=1}^{p-1} a_i\right)^* (a_p a_p^*)$$
$$= \left(\prod_{i=1}^{p-1} (a_i a_i^*)\right)(a_p a_p^*) = \prod_{i=1}^{p} (a_i a_i^*)$$

Следовательно утверждение теоремы верно для $m = p$. □

## 2. Поле частных алгебры скаляров

Пусть $D$ - коммутативное кольцо. Пусть $A$ - $D$-алгебра с сопряжением. Согласно определению [7]-4.2, алгебра скаляров $\operatorname{Re} A$ является коммутативным ассоциативным кольцом. Пусть кольцо $\operatorname{Re} A$ является целостным.[2] Тогда существует поле $F$ частных кольца $\operatorname{Re} A$.[3]

Согласно построениям, выполненным в разделах [8]-4.4.2, [8]-4.4.3, диаграмма представлений $\operatorname{Re} A$-алгебры $A$ имеет вид

$$\operatorname{Re} A \xrightarrow{f_{1,2}}_* A \xrightarrow{f_{2,3}}_* A \qquad f_{1,2}(d): a \to d\,a$$
$$\uparrow_{*f_{1,2}} \qquad f_{2,3}(a): b \to C_A(a,b)$$
$$\operatorname{Re} A \qquad C_A \in \mathcal{L}(A^2; A)$$

Диаграмма представлений $F$-алгебры $G$ имеет вид

$$F \xrightarrow{g_{1,2}}_* B \xrightarrow{g_{2,3}}_* B \qquad g_{1,2}(d): a \to d\,a$$
$$\uparrow_{*g_{1,2}} \qquad g_{2,3}(a): b \to C_B(a,b)$$
$$F \qquad C_B \in \mathcal{L}(B^2; B)$$

Мы определим $F$-алгебру $B$ так, что существует линейный гомоморфизм[4] $\operatorname{Re} A$-алгебры $A$ в $F$-алгебру $B$

$$r_1 : \operatorname{Re} A \to F \quad r_2 : A \to B$$

такой, что гомоморфизм колец $r_1$ является вложением кольца $\operatorname{Re} A$ в кольцо $F$ (страница [1]-86)

(2.1) $$r_1(d) = d/1$$

и образ базиса $\overline{\overline{e}}_A$ $\operatorname{Re} A$-модуля $A$ при отображении $r_2$ является базисом $F$-векторного пространства $B$

(2.2) $$r_2 \circ \overline{e}_{A \cdot \underline{i}} = \overline{e}_{B \cdot \underline{i}}$$

---

[2]Смотри определение целостного кольца на странице [1]-79.
[3]Построение поля частных рассмотренно в [1], страницы 85 - 88.
[4]Смотри определение [6]-6.2.

На основе равенства (2.1) мы можем отождествить $d \in \operatorname{Re} A$ и его образ

$$r_1(d) = d$$

**Теорема 2.1.** *$F$-алгебра $B$ является алгеброй с сопряжением. Поле $F$ является алгеброй скаляров $F$-алгебры $B$. Структурные константы $F$-алгебра $B$ совпадают со структурными константами $\operatorname{Re} A$-алгебры $A$*

$$C_{B\cdot ij}^{\phantom{B\cdot}k} = C_{A\cdot ij}^{\phantom{A\cdot}k} \tag{2.3}$$

*Доказательство.* Из равенства (2.2) следует

$$r_{2\cdot i}^{\phantom{2\cdot}k} = \delta_i^k \tag{2.4}$$

Из равенства (2.4) и теоремы [6]-6.4 следует

$$C_{A\cdot ij}^{\phantom{A\cdot}l} = r_1(C_{A\cdot ij}^{\phantom{A\cdot}k})\delta_k^l = \delta_i^p \delta_j^q C_{B\cdot pq}^{\phantom{B\cdot}l} \tag{2.5}$$

Равенство (2.2) следует из равенства (2.5). Из равенства (2.5) и теоремы [7]-4.5 следует, что $F$-алгебра $B$ является алгеброй с сопряжением и поле $F$ является алгеброй скаляров алгеброй скаляров алгебры $B$. $\square$

В дальнейшем мы будем полагать, что $\operatorname{Re} A$ является полем.

Согласно теореме [7]-4.9

$$aa^* \in \operatorname{Re} A \tag{2.6}$$

В отличие от поля комплексных чисел и алгебры кватернионов, поле $\operatorname{Re} A$ может быть отлично от поля действительных чисел. В поле $\operatorname{Re} A$ может отсутствовать понятие порядка. Поэтому мы не можем интерпретировать выражение (2.6) как норму в алгебре $A$. Более того, это выражение может быть равно 0.

**Теорема 2.2.** *$a \in A$ обратим в $\operatorname{Re} A$-алгебре $A$ тогда и только тогда, когда*

$$aa^* \neq 0 \tag{2.7}$$

*Если условие (2.6) выполнено, то*

$$a^{-1} = \frac{1}{aa^*} a^* \tag{2.8}$$

*Доказательство.* По определению, $a \in A$ обратим, если существует $a^{-1} \in A$ такой, что

$$aa^{-1} = 1 \tag{2.9}$$

Из равенства

$$aa^* = aa^*$$

и равенства (2.6) следует

$$\frac{1}{aa^*} aa^* = \frac{1}{aa^*} aa^* \tag{2.10}$$

Мы можем представить правую часть равенства (2.10) в виде

$$\frac{1}{aa^*} aa^* = \frac{1}{aa^*} \frac{aa^*}{1} = 1 \tag{2.11}$$

Поскольку произведение в $\operatorname{Re} A$-алгебре $A$ - билинейное отображение, то мы можем представить левую часть равенства (2.10) в виде

$$\frac{1}{aa^*} aa^* = a\left(\frac{1}{aa^*} a^*\right) \tag{2.12}$$

Из равенств (2.10), (2.11), (2.12), следует

$$(2.13) \qquad a\left(\frac{1}{aa^*}a^*\right) = 1$$

(2.8) следует из равенств (2.9), (2.13).

Если

$$(2.14) \qquad aa^* = 0$$

то $a \in A$ является левым и правым делителем нуля.[5] Согласно теореме [3]-6.3, $a \in A$ не имеет обратного. □

*Замечание* 2.3. Используя обозначения, рассмотренные в начале раздела, мы можем утверждать, что из теоремы 2.2 следует, что $F$-алгебра $B$ является алгеброй, в которой обратимо максимально возможное множество элементов $\operatorname{Re} A$-алгебры $A$. $F$-алгебра $B$ называется **алгеброй частных алгебры с сопряжением** $A$. □

**Определение 2.4.** Обозначим

$$A_0 = \{a \in A : a\,a^* = 0\}$$

**множество нулей алгебры** $A$. Согласно теореме 2.2, $a \in A_0$ тогда и только тогда, когда либо $a = 0$, либо $a$ является делителем нуля. Обозначим

$$A_1 = \{a \in A : a\,a^* \neq 0\}$$

**множество обратимых элементов алгебры** $A$. □

Пусть $a \in A_1$, $b \in A$. **Левая дробь** представлена выражением

$$a^{-1}b = \frac{1}{a\,a^*}a^*b$$

**Правая дробь** представлена выражением

$$ba^{-1} = \frac{1}{a\,a^*}b\,a^*$$

Множество дробей в алгебре не ограничено левыми или правыми дробями. Например, выражения

$$(a^{-1}b)(c^{-1}d) \quad a^{-1}b + c^{-1}d$$

также являются дробями.

На множестве дробей можно определить несколько отношений эквивалентности. Например, если $d \in \operatorname{Re} A$, то

$$a^{-1}b = (da)^{-1}(db) \quad (a^{-1}(db))(c^{-1}f) = (a^{-1}b)(c^{-1}(df))$$

Однако вопрос о канонической форме дроби - задача нетривиальная, по крайней мере, в данный момент времени.

---

[5]Из равенств

$$aa^* = C^{\mathbf{0}}_{\mathbf{00}}a^{\mathbf{0}}a^{\mathbf{0}} - C^{\mathbf{0}}_{kl}a^k a^l$$
$$a^*a = C^{\mathbf{0}}_{\mathbf{00}}a^{\mathbf{0}}a^{\mathbf{0}} - C^{\mathbf{0}}_{kl}a^k a^l$$

и равенства (2.14) следует $aa^* = 0$.

## 3. Алгебра многочленов

Пусть $D$ - коммутативное кольцо характеристики $0$. Пусть $A$ - $D$-алгебра. **Алгебра многочленов** $A[x]$ **над** $D$**-алгеброй** $A$ порождена множеством одночленов. Структура одночлена $p_k$ степени $k$, $k > 0$, одной переменной над ассоциативной $D$-алгеброй $A$ описана в следующей теореме (раздел [5]-5.2).

**Теорема 3.1.** *Одночлен степени $0$ имеет вид $a_0$, $a_0 \in A$. Для $k > 0$,*

$$p_k(x) = p_{k-1}(x) x a_k$$

*где $a_k \in A$.*

*Доказательство.* Действительно, последний множитель одночлена $p_k(x)$ является либо $a_k \in A$, либо имеет вид $x^l$, $l \geq 1$. В последнем случае мы положим $a_k = 1$. Множитель, предшествующий $a_k$, имеет вид $x^l$, $l \geq 1$. Мы можем представить этот множитель в виде $x^{l-1}x$. Следовательно, утверждение доказано. □

В частности, одночлен степени $1$ имеет вид $p_1(x) = a_0 x a_1$. Из теорем 3.1, [5]-3.41 следует, что каждому одночлену $p_k$ мы можем сопоставить тензор $a_0 \otimes a_1 \otimes ... \otimes a_k$.

В неассоциативной алгебре порядок сомножителей становится существенным. Поэтому теорема 3.1 приобретает следующую форму.

**Теорема 3.2.** *Одночлен степени $0$ имеет вид $a_0$, $a_0 \in A$. Для $k > 0$, существуют одночлены $p_l$, $p_m$, $m + l = k$, такие, что*

$$p_k(x) = a_{k \cdot 1} p_l(x) a_{k \cdot 2} p_m(x) a_{k \cdot 3}$$

*где $a_{k \cdot 1}$, $a_{k \cdot 2}$, $a_{k \cdot 3} \in A$.* □

Если $D$-алгебра $A$ является алгеброй с сопряжением, то мы можем распространить отображение сопряжения на алгебру многочленов, а также рассматривать многочлены над кольцом $\operatorname{Re} A$. Так как структура многочлена над кольцом $\operatorname{Re} A$ отличается от структуры многочлена над алгеброй $A$, то определение связи между алгебрами $\operatorname{Re} A[x]$ и $A[x]$ является нетривиальной задачей.

Если $p(x)$ - одночлен над алгеброй, то, согласно теореме 1.4, мы можем рассматривать выражение $p(x)(p(x))^*$ как многочлен $r(y)$ с переменной $y = xx^*$ над алгеброй $\operatorname{Re} A$. Хотя для произвольного многочлена $p(x) \in A[x]$ выражения $p(x) + (p(x))^*$, $p(x)(p(x))^*$ принимают значения в кольце $\operatorname{Re} A$, совсем не очевидно, можем ли мы эти выражения рассматривать как многочлены над кольцом $\operatorname{Re} A$.

Существуют алгебры с сопряжением, в которых сопряжение не зависит линейно от тождественного отображения (смотри, например, раздел [2]-6). В этом случае для произвольного полинома выражение

$$p(x)(p(x))^*$$

зависит от двух переменных: $x$ и $x^*$.

Рассмотрим алгебры с сопряжением, в которых сопряжение линейно зависит от тождественного отображения

$$x^* = s(x) = s_{i \cdot 0} \, x \, s_{i \cdot 1}$$

(например, отображение [4]-(4.3.35) в алгебре кватернионов, отображение [4]-(4.5.99) в алгебре октонионов). В этом случае для произвольного полинома $p(x)$ выражение

$$p(x)(p(x))^* = r(x)$$

является полиномом.

Пусть кольцо $\operatorname{Re} A$ является полем.

**Определение 3.3.** Пусть $p(x)$ - полином над алгеброй $A$. $a \in A$ называется **корнем полинома** $p$, если $p(a) \in A_0$. □

Согласно теореме 2.2, для полинома $p(x)$ определено выражение

$$\text{(3.1)} \qquad \frac{1}{p(x)(p(x))^*} \in \operatorname{Re} A$$

для любого $x \in A$, отличного от корня полинома $p(x)$. Следовательно, отображение

$$\text{(3.2)} \qquad (p(x))^{-1} = \frac{1}{p(x)(p(x))^*}(p(x))^*$$

корректно определено для $x$, не являющихся корнями многочлена $p$.

Алгебра $A(x)$, порождённая выражениями вида (3.2) называется **алгеброй рациональных отображений алгебры** $A$.

## 4. Идеал

**Определение 4.1.** Подгруппа $B$ аддитивной группы алгебры с сопряжением $A$ называется **левым идеалом алгебры**,[6] если

$$aB \subset B \quad a \in A$$

Подгруппа $B$ аддитивной группы алгебры с сопряжением $A$ называется **правым идеалом алгебры**, если

$$Ba \subset B \quad a \in A$$

Подгруппа $B$ аддитивной группы алгебры с сопряжением $A$ называется **идеалом алгебры**, если $B$ одновременно является левым и правым идеалом. □

**Пример 4.2.** Пусть $A$ - алгебра с сопряжением, $a \in A$. Множество $Aa$ является левым идеалом, называемым **левым главным идеалом** алгебры $A$.[7] Множество $aA$ является правым идеалом, называемым **правым главным идеалом** алгебры $A$. Множество $AaA$ является идеалом, называемым **главным идеалом** алгебры $A$. □

**Теорема 4.3.** *Пусть $A$ - ассоциативная алгебра с сопряжением, $a \in A_0$. Тогда*

$$\text{(4.1)} \qquad Aa \in A_0$$

$$\text{(4.2)} \qquad aA \in A_0$$

---

[6]Это определение опирается на определение [1], стр. 75.

[7][1], страница 75.

*Доказательство.* Пусть $b \in A$. Из определения 2.4 следует

(4.3) $\quad (ba)(ba)^* = (ba)(a^*b^*) = b(aa^*)b^* = 0$

(4.4) $\quad (ab)(ab)^* = (ab)(b^*a^*) = a(bb^*)a^* = (bb^*)(aa^*) = 0$

Утверждение (4.1) следует из равенства (4.3). Утверждение (4.2) следует из равенства (4.4). $\square$

**Теорема 4.4.** *Пусть $A$ - ассоциативная алгебра с сопряжением, $a \in A$, $b \in A_0$. Алгебра многочленов $A[x]$ имеет левый идеал*

$$Z^1_l(a,b)A[x] = \{p \in A[x] : p(a) \in Ab\}$$

*Алгебра многочленов $A[x]$ имеет правый идеал*

$$Z^1_r(a,b)A[x] = \{p \in A[x] : p(a) \in bA\}$$

*Доказательство.* Теорема является следствием определений, рассмотренных в примере 4.2 и теоремы 4.3. $\square$

Верны также аналогичные теоремы.

**Теорема 4.5.** *Пусть $A$ - ассоциативная алгебра с сопряжением, $a \in A$. Алгебра многочленов $A[x]$ имеет идеал*

$$Z^1(a)A[x] = \{p \in A[x] : p(a) = 0\}$$

$\square$

**Теорема 4.6.** *Пусть $A$ - ассоциативная алгебра с сопряжением, $a \in A$, $b \in A_0$. Алгебра рациональных отображений $A(x)$ имеет левый идеал*

$$Z^1_l(a,b)A(x) = \{p \in A(x) : p(a) \in Ab\}$$

*Алгебра рациональных отображений $A(x)$ имеет правый идеал*

$$Z^1_r(a,b)A(x) = \{p \in A(x) : p(a) \in bA\}$$

$\square$

**Теорема 4.7.** *Пусть $A$ - ассоциативная алгебра с сопряжением, $a \in A$. Алгебра рациональных отображений $A(x)$ имеет идеал*

$$Z^1(a)A(x) = \{p \in A(x) : p(a) = 0\}$$

$\square$

## 5. Список литературы

## 6. Предметный указатель

# 7. Специальные символы и обозначения